\documentclass[amstex,12pt, amssymb]{article}

\usepackage{mathtext}
\usepackage[cp1251]{inputenc}
\usepackage[T2A]{fontenc}
\usepackage[dvips]{graphicx}
\usepackage{amsmath}
\usepackage{amssymb}
\usepackage{amsxtra}
\usepackage{latexsym}
\usepackage{ifthen}

\newcounter{lemma}[section]

\newcounter{corollary}[section]

\newcounter{remark}[section]

\newcounter{theorem}[section]

\newcounter{proposition}[section]

\newcounter{example}

\numberwithin{equation}{section}

\begin{document}

\markboth{\centerline{O.~DOVHOPIATYI, E.~SEVOST'YANOV}}
{\centerline{ON SOLUTIONS OF THE BELTRAMI EQUATIONS...}}

\def\cc{\setcounter{equation}{0}
\setcounter{figure}{0}\setcounter{table}{0}}

\overfullrule=0pt


\author{OLEKSANDR DOVHOPIATYI, EVGENY SEVOST'YANOV\\}

\title{
{\bf ON BELTRAMI EQUATIONS WITH INVERSE CONDITIONS AND HYDRODYNAMIC
NORMALIZATION}}

\date{\today}
\maketitle

\begin{abstract}
We consider problems concerning the existence of solutions of the
Beltrami equations and their convergence in the entire complex
plane. We are mainly interested in the case when these solutions
satisfy the so-called hydrodynamic normalization condition in the
neighborhood of infinity. Under conditions related to dilations of
inverse mappings, we have established the existence of such
solutions in the class of continuous Sobolev mappings. We have also
obtained results on the locally uniform limit of a sequence of such
solutions.
\end{abstract}


\section{Introduction}

In the publications of the second co-author, some problems related
to the existence of continuous Sobolev solutions of the Beltrami
equations were considered (see \cite{Sev}, \cite{SevSkv}). Note
that, as a rule, we were talking about solutions defined in the unit
disk and satisfying the normalization $f(0)=0$ and $f(1)=1.$ In this
manuscript, we will consider the problem of the existence of
solutions to the Beltrami equations with a different normalization.
Some publications that are close in context to the problem under
study should be mentioned. In particular, existence theorems for
solutions of homeomorphic solutions of the Beltrami equations with
hydrodynamic normalization were obtained in~\cite{GRSY}, where
slightly different conditions on the dilatation of the mapping were
used. Also, in \cite{DS$_1$}--\cite{DS$_2$} we have studied the
problem of compactness of solutions of the Beltrami equations with
hydrodynamic normalization.

\medskip
Let $D$ be a domain in ${\Bbb C}.$ In what follows, a mapping
$f:D\rightarrow{\Bbb C}$ is assumed to be {\it sense-preserving,}
moreover, we assume that $f$ has partial derivatives almost
everywhere. Put $f_{\overline{z}} = \left(f_x + if_y\right)/2$ and
$f_z = \left(f_x - if_y\right)/2.$ Let $\mu=\mu(z):D \rightarrow
{\Bbb D}$ and $\nu=\nu(z):D \rightarrow {\Bbb D}$ be Lebesgue
measurable functions. We set
\begin{equation}\label{eq1}
K_{\mu, \nu}(z)=\quad\frac{1+|\mu
(z)|+|\nu(z)|}{1-|\mu\,(z)|-|\nu(z)|}\,.
\end{equation}
The {\it complex dilatation} of $f$ at $z\in D$ is defined as
follows: $\mu(z)=\mu_f(z)=f_{\overline{z}}/f_z$ for $f_z\ne 0$ and
$\mu(z)=0$ otherwise. The {\it maximal dilatation} of $f$ at $z$ is
the following function:
\begin{equation}\label{eq1D}
K_{\mu}(z)=K_{\mu_f}(z)=\quad\frac{1+|\mu (z)|}{1-|\mu\,(z)|}\,.
\end{equation}
Note that the Jacobian of $f$ at $z\in D$ may be calculated
according to the relation
$$J(z,
f)=|f_z|^2-|f_{\overline{z}}|^2\,.$$ Since we assume that the map
$f$ is sense preserving, the Jacobian of this map is positive at all
points where $f$ is differentiable. Let ${\Bbb D}=\{z\in {\Bbb C}:
|z|<1\},$ and let $\mu:D\rightarrow {\Bbb D}$ be a Lebesgue
measurable function. We define the {\it maximal dilatation}
corresponding to a complex dilatation~$\mu$ by~(\ref{eq1D}).
It is easy to see that
$$K_{\mu_f}(z)=\frac{|f_z|+|f_{\overline{z}}|}{|f_z|-|f_{\overline{z}}|}$$
whenever partial derivatives of $f$ exist at $z\in D$ and, in
addition, $J(z, f)\ne 0.$
Set $\Vert f^{\,\prime}(z)\Vert=|f_z|+|f_{\overline{z}}|.$ Recall
that a homeomorphism $f:D\rightarrow {\Bbb C}$ is said to be {\it
quasiconformal} if $f\in W_{\rm loc}^{1, 2}(D)$ and, in addition,
$\Vert f^{\,\prime}(z)\Vert^2\leqslant K\cdot |J(z, f)|$ for some
constant $K\geqslant 1$ almost everywhere.

\medskip
A {\it Beltrami equation with two characteristics} is a differential
equation of the form
\begin{equation}\label{eq2J}
f_{\overline{z}}=\mu(z)\cdot f_z+\nu(z)\cdot \overline{f_z}\,,
\end{equation}
where $\mu=\mu(z)$ and $\nu=\nu(z)$ are given measurable functions
with $|\mu(z)|<1$ and $|\nu(z)|<1$  a.a. Let $\mu=\mu(z):
D\rightarrow {\Bbb D}$ and $\nu=\nu(z): D\rightarrow {\Bbb D}$ be
functions such that the relation $|\mu(z)|+|\nu(z)|< 1$ holds for
almost any $z\in D.$ We will consider that $\mu(z)=\nu(z)\equiv 0$
for any $z\in {\Bbb C}\setminus D.$ Fix $n\geqslant 1$ and set
\begin{equation}\label{eq12:} \mu_n(z)= \left
\{\begin{array}{rr}
 \mu(z),&  z\in {\Bbb C}, K_{\mu, \nu}(z)\leqslant n,
\\ 0\ , & \text{otherwise\,in \,} {\Bbb C}\,,
\end{array} \right.
\end{equation}
and
\begin{equation}\label{eq12:A} \nu_n(z)= \left
\{\begin{array}{rr}
 \nu(z),&  z\in {\Bbb C}, K_{\mu, \nu}(z)\leqslant n,
\\ 0\ , & \text{otherwise\,in \,} {\Bbb C}\,.
\end{array} \right.
\end{equation}
Let $f_n:{\Bbb C}\rightarrow {\Bbb C}$ be a homeomorphic solution of
the equation $(f_n)_{\overline{z}}=\mu_n(z)\cdot
(f_n)_z+\nu_n(z)\cdot \overline{(f_n)_z}$ (it exists by
\cite[Theorem~9.2]{GRSY}). Set $g_n(z):=f^{\,-1}_n(z).$ Observe that
$g_n:{\Bbb C}\rightarrow {\Bbb C}$ is quasiconformal, in particular,
$g_n$ is almost everywhere differentiable in ${\Bbb C}.$ Observe
that $f_n$ is conformal at the neighborhood of the infinity, so
there is a continuous extension $f:\overline{\Bbb C}\rightarrow
\overline{\Bbb C}.$ Thus $f_n({\Bbb C})={\Bbb C}$ and
$f_n(\infty)=\infty.$ By Proposition~2.1 in \cite{GRSY},
$f_n(z)=a_nz+b_n+o(1)$ as $z\rightarrow\infty,$ where $a_n, b_n\in
{\Bbb C}$ and $a_n\ne 0.$ We may consider that $a_n=1$ and $b_n=0$
for any $n\in{\Bbb N},$ because
$\widetilde{f}_n(z):=\frac{1}{a_n}\cdot f_n(z)-b_n$ satisfies the
same equation $f_{\overline{z}}=\mu(z)\cdot f_z+\nu(z)\cdot
\overline{f_z}$ in ${\Bbb C}.$

\medskip
Note that such a function $\mu$ is unique. Indeed, the same function
is a solution to the usual Beltrami equation
$f_{\overline{z}}=\mu^*_n(z)\cdot f_z,$ where
$\mu_1(z)=\mu_n(z)+\nu_n(z)\cdot \frac{\overline{f_z}}{f_z},$ and,
by the triangle inequality, $|\mu^*_n(z)|\leqslant
|\mu_n(z)|+|\nu_n(z)|\leqslant \frac{n-1}{n+1}<1.$ Thus, $f$ is
unique by~\cite[Theorem~20.4.15]{AIM}.

\medskip
Let $K_{\mu_{g_n}}(w)$ be a maximal dilatation of $g_n,$ namely,
\begin{equation}\label{eq3D}
K_{\mu_{g_n}}(w)=\quad
\frac{{|(g_n)_w|}^2-
{|(g_n)_{\overline{w}}|}^2}{{(|(g_n)_w|-|(g_n)_{\overline{w}}|)}^2}\,.
\end{equation}
We also define the {\it inner dilatation of $g_n$ of the order $p$
at a point $w$} by the equation
\begin{equation}\label{eq18}
K_{I, p}(w,
g_n)=\frac{{|(g_n)_w|}^2-{|(g_n)_{\overline{w}}|}^2}{{(|(g_n)_w|-|(g_n)_{\overline{w}}|)}^p}\,.
\end{equation}

\medskip
The following statement holds.

\medskip
\begin{theorem}\label{th1A}{\sl\,
Let $D$ be a domain in ${\Bbb C}$ such that $\overline{D}$ is a
compact set in ${\Bbb C},$ let $\mu=\mu(z):{\Bbb C}\rightarrow {\Bbb
D}$ and $\nu=\nu(z):{\Bbb C}\rightarrow {\Bbb D}$ be Lebesgue
measurable functions vanishing outside $D$ such that the relation
$|\mu(z)|+|\nu(z)|< 1$ holds for almost any $z\in D.$ In addition,
let $\mu_n,$ $\nu_n,$ $f_n$ and $g_n$ as above, $n=1,2,\ldots, .$
Let $Q:{\Bbb C}\rightarrow[1, \infty]$ be a Lebesgue measurable
function. Assume that the following conditions hold:

\medskip
1) for each $0<r_1<r_2<1$ and $y_0\in {\Bbb C}$  there is a set
$E\subset[r_1, r_2]$ of positive linear Lebesgue measure such that
the function $Q$ is integrable over the circles $S(y_0, r)$ for any
$r\in E;$

\medskip
2) there exist a number $1<p\leqslant 2$ such that, for any bonded
domain $G\subset {\Bbb C}$ there exists a constant $M=M_G>0$ such
that
\begin{equation}\label{eq10B}
\int\limits_{G}K_{I, p}(w, g_n)\,dm(w)\leqslant M
\end{equation}
for all $n=1,2,\ldots ,$ where $K_{I, p}(w, g_n)$ is defined
in~(\ref{eq18});

\medskip
3) the inequality
\begin{equation}\label{eq10C}
K_{\mu_{g_n}}(w)\leqslant Q(w)
\end{equation}
holds for a.e. $w\in {\Bbb C},$ where $K_{\mu_{g_n}}$ is defined
in~(\ref{eq3D}).
Then the equation~(\ref{eq2J}) has a continuous $W_{\rm loc}^{1,
p}({\Bbb C})$-solution $f$ in ${\Bbb C}$ such that
$f(z)=z+\varepsilon(z),$ where $\varepsilon(z)\rightarrow 0$ as
$z\rightarrow\infty.$ }
\end{theorem}

\medskip
\begin{corollary}\label{cor5}
{\sl\, In particular, the conclusion of Theorem~\ref{th1A} holds if,
in this theorem, we abandon condition~1), accept condition~3), and
replace condition~2) with the requirement $Q\in L_{\rm loc}^1({\Bbb
C}).$ If $G$ is some bounded domain in ${\Bbb C}$ and $K$ is a
compactum in $G,$ then there is some domain $G^{\,\prime}\subset
{\Bbb C}$ and a function $Q^{\,\prime}$ which equal to $Q$ in
$G^{\,\prime}$ and vanishing outside $G^{\,\prime}$ such that
$Q^{\,\prime}$ is integrable in ${\Bbb C}$ and the relation
\begin{equation}\label{eq2E}
|f(x)-f(y)|\leqslant\frac{C}{\log^{1/2}\left(1+\frac{r_0}{2|x-y|}\right)}
\end{equation}
holds for any $x, y\in K,$ where $C=C(K, \Vert Q^{\,\prime}\Vert_1,
G)>0$ is some constant depending only on $K,$ $G$ and $\Vert
Q^{\,\prime}\Vert_1,$ $\Vert Q^{\,\prime}\Vert_1$ denotes $L^1$-norm
of $Q^{\,\prime}$ in ${\Bbb R}^n,$ and $r_0=d(K,
\partial G).$
}
\end{corollary}

\section{Proof of the main results}

Sometimes, instead of ${\rm dist}\,(A, B),$ we also write $d(A, B),$
if a misunderstanding is impossible. Given sets $E$ and $F$ and a
given domain $D$ in $\overline{{\Bbb R}^n}={\Bbb R}^n\cup
\{\infty\},$ we denote by $\Gamma(E, F, D)$ the family of all paths
$\gamma:[0, 1]\rightarrow \overline{{\Bbb R}^n}$ joining $E$ and $F$
in $D,$ that is, $\gamma(0)\in E,$ $\gamma(1)\in F$ and
$\gamma(t)\in D$ for all $t\in (0, 1).$ Everywhere below, unless
otherwise stated, the boundary and the closure of a set are
understood in the sense of the extended Euclidean space
$\overline{{\Bbb R}^n}.$ Let $x_0\in\overline{D},$ $x_0\ne\infty,$
\begin{equation}\label{eq1E}
B(x_0, r)=\{x\in {\Bbb R}^n: |x-x_0|<r\}\,, \quad {\Bbb B}^n=B(0,
1)\,, \end{equation}
$$S(x_0,r) = \{
x\,\in\,{\Bbb R}^n : |x-x_0|=r\}\,, S_i=S(x_0, r_i)\,,\quad
i=1,2\,,$$
\begin{equation*}\label{eq1**}
A=A(x_0, r_1, r_2)=\{ x\,\in\,{\Bbb R}^n : r_1<|x-x_0|<r_2\}\,.
\end{equation*}
Let $f:D\rightarrow{\Bbb R}^n,$ $n\geqslant 2,$ and let $Q:{\Bbb
R}^n\rightarrow [0, \infty]$ be a Lebesgue measurable function such
that $Q(y)\equiv 0$ for $y\in{\Bbb R}^n\setminus f(D).$ Let
$A=A(y_0, r_1, r_2).$ Let $\Gamma_f(y_0, r_1, r_2)$ denotes the
family of all paths $\gamma:[a, b]\rightarrow D$ such that
$f(\gamma)\in \Gamma(S(y_0, r_1), S(y_0, r_2), A(y_0, r_1, r_2)),$
i.e., $f(\gamma(a))\in S(y_0, r_1),$ $f(\gamma(b))\in S(y_0, r_2),$
and $f(\gamma(t))\in A(y_0, r_1, r_2)$ for any $a<t<b.$ We say that
{\it $f$ satisfies the inverse Poletsky inequality} at $y_0\in f
(D)$ if the relation
\begin{equation}\label{eq2*A}
M(\Gamma_f(y_0, r_1, r_2))\leqslant \int\limits_A Q(y)\cdot \eta^n
(|y-y_0|)\, dm(y)
\end{equation}
holds for any $0<r_1<r_2<r_0:=\sup\limits_{y\in
D^{\,\prime}}|y-y_0|$ and any Lebesgue measurable function
$\eta:(r_1, r_2)\rightarrow [0, \infty]$ such that
\begin{equation}\label{eqA2}
\int\limits_{r_1}^{r_2}\eta(r)\,dr\geqslant 1 \,.
\end{equation}

\medskip
For domains $D, D^{\,\prime}\subset {\Bbb R}^n,$ $n\geqslant 2,$ and
a Lebesgue measurable function $Q:{\Bbb R}^n\rightarrow [0, \infty]$
equal to zero outside the domain $D^{\,\prime},$ we define by
$\frak{R}_Q(D, D^{\,\prime})$ the family of all open discrete
mappings $f:D\rightarrow D^{\,\prime}$ such that
relation~(\ref{eq2*A}) holds for each point $y_0\in D^{\,\prime}.$
Note that the definition of class $\frak{R}_Q(D, D^{\,\prime})$ does
not require that the domain $D$ be mapped onto the domain
$D^{\,\prime}$ surjectively under the mapping $f\in \frak{R}_Q(D,
D^{\,\prime}).$ In what follows, $\mathcal{H}^{n-1}$ denotes
$(n-1)$-dimensional Hausdorff measure. The following statement holds
(see~\cite[Theorem~1.1]{SevSkv}).

\medskip
\begin{proposition}\label{pr1}
{\sl Let $D$ and $D^{\,\prime}$ be domains in ${\Bbb R}^n,$
$n\geqslant 2,$ and let $D^{\,\prime}$ be a bounded domain. Suppose
that, for each point $y_0\in D^{\,\prime}$ and for every
$0<r_1<r_2<r_0:=\sup\limits_{y\in D^{\,\prime}}|y-y_0|$ there is a
set $E\subset[r_1, r_2]$ of a positive linear Lebesgue measure such
that the function $Q$ is integrable with respect to
$\mathcal{H}^{n-1}$ over the spheres $S(y_0, r)$ for every $r\in E.$
Then the family of mappings $\frak{R}_Q(D, D^{\,\prime})$ is
equicontinuous at each point $x_0\in D.$ }
\end{proposition}

\medskip
\begin{remark}\label{rem1}
For domains $D, D^{\,\prime}\subset {\Bbb R}^n,$ $n\geqslant 2,$ and
a Lebesgue measurable function $Q:{\Bbb R}^n\rightarrow [0, \infty]$
equal to zero outside the domain $D^{\,\prime},$ we define by
$\frak{R}_Q(D, D^{\,\prime})$ the family of all open discrete
mappings $f:D\rightarrow D^{\,\prime}$ such that relation
$$
M(\Gamma_f(y_0, r_1, r_2))\leqslant \frac{\omega_{n-1}}{I^{n-1}}$$
holds for each point $y_0\in D^{\,\prime}$ with
$I=\int\limits_{r_1}^{r_2}\frac{dt}{tq_{y_0}^{1/(n-1)}(t)},$
$q_{y_0}(r)=\frac{1}{\omega_{n-1}r^{n-1}}\int\limits_{S(y_0,
r)}Q(y)\,d\mathcal{H}^{n-1}(y),$  $\omega_{n-1}$ is the area of the
unit sphere ${\Bbb S}^{n-1}$ in ${\Bbb R}^n,$  and any
$0<r_1<r_2<r_0:=\sup\limits_{y\in D^{\,\prime}}|y-y_0|.$ {\it
Suppose that, for each point $y_0\in D^{\,\prime}$ and for every
$0<r_1<r_2<r_0:=\sup\limits_{y\in D^{\,\prime}}|y-y_0|$ there is a
set $E\subset[r_1, r_2]$ of a positive linear Lebesgue measure such
that the function $Q$ is integrable with respect to
$\mathcal{H}^{n-1}$ over the spheres $S(y_0, r)$ for every $r\in E.$
Then the family of mappings $\frak{R}_Q(D, D^{\,\prime})$ is
equicontinuous at each point $x_0\in D.$}

\medskip
The proof of this assertion repeats almost exactly the proof of
Theorem 1.1 in~\cite{SevSkv} and is therefore omitted. Note that the
function $Q$ in the above statement can be extended outside the
domain $D^{\,\prime}$ arbitrarily, and not necessarily by zero.
\end{remark}

\medskip
For the case of mappings of a domain $D$ onto itself, the following
lemma is proved in~\cite[Lemma~2]{DS$_3$},
cf.~\cite[Theorem~9.1]{GRSY} and \cite[Lemma~5.1]{SevSkv}. We need
to formulate it in a somewhat more general case, when the mappings
do not, in general, map the domain $D$ onto itself. We note that the
proof of this assertion is quite similar to the proof of the above
particular case, however, we present it in full in the text.

\medskip
\begin{lemma}\label{pr2}
{\sl\, Let $1<p\leqslant 2,$ and let $f_n,$ $n=1,2,\ldots $ -- be a
sequence of sense-preserving homeomorphisms of a domain $D$ into
${\Bbb R}^n$ which belong to the class $W_{\rm loc}^{1, 2}(D)$ and
satisfying the equation
\begin{equation}\label{eq1B}
\overline{\partial}f_n=\partial f_n\mu_n(z)+\overline{\partial
f_n}\nu_n(z)\,,
\end{equation}
where $\mu_n$ and $\nu_n$ are Lebesgue measurable functions
satisfying the relation $|\nu_n(z)|+|\mu_n(z)|<1$ almost everywhere.
Assume that, $f_n$ converge to a mapping $f:D\rightarrow {\Bbb C}$
locally uniformly as $n\rightarrow \infty,$ and $\mu_n(z)$ and
$\nu_n(z)$ converge to $\mu(z)$ and $\nu(z)$ as $n\rightarrow\infty$
a.e.  Assume that, the inverse mappings $g_n:=f_n^{\,-1}$ belong to
$W_{\rm loc}^{1, 2}(f_n(D)),$ while
$$\int\limits_{f_n(D)}K_{I, p}(w, g_n)\,dm(w)\leqslant M$$
for some $M>0$ and any $n=1,2,\ldots .$
Then $f\in W_{\rm loc}^{1, p}(D)$ and, in addition, $\mu$ and $\nu$
are complex characteristics of the mapping $f,$ in other words,
$\overline{\partial}f=\partial f\mu(z)+\overline{\partial f}\nu(z)$
for almost any $z\in D.$
 }
\end{lemma}

\medskip
\begin{proof}
We will generally follow the scheme outlined
in~\cite[Theorem~9.1]{GRSY}, cf.~\cite[Lemma~5.1]{SevSkv}. Denote
$\partial f=f_z$ and $\overline{\partial}f=f_{\overline{z}}.$ Let
$C$ be arbitrary compactum in $D.$ Since $g_n=f_n^{\,-1}\in W_{\rm
loc}^{1, 2},$ by the assumption, $g_n$ have $N$-property of Luzin,
see e.g. \cite[Corollary~B]{MM}. Then $J(z, f_n)$ almost everywhere
is not equal to zero, see, for example, \cite[Theorem~1]{Pon}.
Moreover, the change of variables formula under the integral holds,
see~\cite[Theorem~3.2.5]{Fe}. In this case, we obtain that
$$\int\limits_{C}{\Vert f^{\,\prime}_n(z)\Vert}^p\,dm(z)=
\int\limits_C \frac{{\Vert f^{\,\prime}_n(z)\Vert}^p}{J(z,
f_n)}\cdot J(z, f_n)\,dm(z)=$$
\begin{equation}\label{eq4D}
=\int\limits_{f_n(C)}K_{I, p}(w, g_n)\, dm(w)\leqslant M<\infty\,.
\end{equation}
It follows from~(\ref{eq4D}) that $f\in W_{\rm loc}^{1, p}$ and,
besides that, $\partial f_n$ and $\overline{\partial} f_n$ converge
weakly in $L_{\rm loc}^1(D)$ to $\partial f$ and
$\overline{\partial} f,$ respectively (see~\cite[Lemma~III.3.5]{Re};
cf.~\cite[Lemma~2.1]{RSY$_2$}).

\medskip
It remains to show that $f$ is a solution of the equation
$f_{\overline{z}}=\mu(z)\cdot f_z+\nu(z)\cdot \overline{f_z}.$ Set
$\zeta(z)=\overline{\partial} f(z)-\mu(z)\partial f(z)-\nu
\overline{\partial f(z)}$ and let us to show that $\zeta(z)=0$
almost everywhere. Let $B$ be arbitrary disk belonging together with
its closure to $D.$ By the triangle inequality
\begin{equation}\label{eq9}
\left|\int\limits_B\zeta(z)\,dm(z)\right|\leqslant
I_1(n)+I_2(n)+I_3(n)\,,
\end{equation}
where
\begin{equation}\label{eq7}
I_1(n)=\left|\int\limits_B(\overline{\partial}f(z)-\overline{\partial}f_n(z))\,dm(z)\right|\,,
\end{equation}
\begin{equation}\label{eq8}
I_2(n)=\left|\int\limits_B(\mu(z)\partial f(z)-\mu_n(z)\partial
f_n(z))\,dm(z)\right|
\end{equation}
and
\begin{equation}\label{eq9C}
I_3(n)=\left|\int\limits_B(\nu(z)\overline{\partial
f(z)}-\nu_n(z)\overline{\partial f_n(z)})\,dm(z)\right|\,.
\end{equation}
Due to the mentioned above, $I_1(n)\rightarrow 0$ as
$n\rightarrow\infty.$ It remains to deal with the expressions
$I_2(n)$ and $I_3(n).$  To do this, note that, by the triangle
inequality,
$I_2(n)\leqslant I^{\,\prime}_2(n)+I^{\,\prime\prime}_2(n),$ where
$$I^{\,\prime}_2(n)=
\left|\int\limits_B\mu(z)(\partial f(z)-\partial
f_n(z))\,dm(z)\right|$$
and
$$I^{\,\prime\prime}_2(n)=
\left|\int\limits_B(\mu(z)-\mu_n(z))\partial
f_n(z)\,dm(z)\right|\,.$$
Due to the weak convergence of $\partial f_n\rightarrow
\partial f$ in $L^1_{\rm loc}(D)$ as $n\rightarrow\infty,$ we
obtain that $I^{\,\prime}_2(n)\rightarrow 0$ as
$n\rightarrow\infty,$ because $\mu\in L^{\infty}(D).$ Moreover,
since the above mapping $\partial f$ is integrable in degree $p>1,$
there is absolute continuity in the integral$\int\limits_E|\partial
f(z)|\,dm(z).$ Besides that, since $\partial f_n\rightarrow
\partial f$ weakly in $L^1_{\rm loc}(D),$ given
$\varepsilon>0$ there is $\delta=\delta(\varepsilon)>0$ such that
$$\int\limits_E|\partial f_n(z)|\,dm(z)\leqslant$$
\begin{equation}\label{eq5}\leqslant
\int\limits_E|\partial f_n(z)-\partial f(z)|\,dm(z)+
\int\limits_E|\partial f(z)|\,dm(z)<\varepsilon\,,
\end{equation}
whenever $m(E)<\delta,$ $E\subset B,$ and numbers $n$ are
sufficiently large.

\medskip
Finally, by Egorov’s theorem (see~\cite[Theorem~III.6.12]{Sa}), for
any $\delta>0$ there is a set $S\subset B$ such that $m(B\setminus
S)<\delta$ and $\mu_n(z)\rightarrow \mu(z)$ uniformly in~$S.$ Then
$|\mu_n(z)-\mu(z)|<\varepsilon$ for all $n\geqslant n_0,$ some
$n_0=n_0(\varepsilon)$ and all $z\in S.$ Due to the
conditions~(\ref{eq5}), (\ref{eq4D}) and by  H\"{o}lder's
inequality, we obtain that
$$I^{\,\prime\prime}_2(n)\leqslant \varepsilon \int\limits_S|\partial
f_n(z)|\,dm(z)+ 2\int\limits_{B\setminus S}|\partial
f_n(z)|\,dm(z)<$$
\begin{equation}\label{eq6}
<\varepsilon\cdot\left\{\left(\int\limits_D K_{I, p}(w, g_n)\,
dm(w)\right)^{1/p}\cdot (m(B))^{(p-1)/p}+2\right\}\leqslant
\end{equation}
$$\leqslant \varepsilon\cdot\left\{M^{1/p}\cdot (m(B))^{(p-1)/p}+2\right\}$$
for the same $n\geqslant n_0.$ Thus,
$I_2^{\,\prime\prime}(n)\rightarrow 0$ as $n\rightarrow\infty,$
therefore $I_2(n)\rightarrow 0$ as $n\rightarrow\infty.$ The fact
that
\begin{equation}\label{eq9A}
I_3(n)\rightarrow 0
\end{equation}
as $n\rightarrow\infty$ may be proved similarly. Thus,
by~(\ref{eq7}), (\ref{eq8}), (\ref{eq9C}), (\ref{eq6}) and
(\ref{eq9A}) it follows that $\int\limits_B\zeta(z)\,dm(z)=0$ for
any disks $B,$ compactly embedded in $D.$ Based on the Lebesgue
theorem on differentiation of an indefinite integral
(see~\cite[IV(6.3)]{Sa}), it follows that $\zeta(z)=0$ almost
everywhere in $D.$ The lemma is proved.~$\Box$
\end{proof}

\medskip
{\it Proof of Theorem~\ref{th1A}}. Let $f_n,$ $n=1,2,\ldots,$ be a
mapping from the condition of the theorem.
Let us to prove that $\{f_n\}^{\infty}_{n=1}$ forms a normal family
of mappings. Fix an arbitrary compact set $C\subset {\Bbb C}.$ Since
$\overline{D}$ is a compactum in ${\Bbb C},$ there is a domain
$G\subset {\Bbb C}$ with a compact closure in ${\Bbb C}$ such that
$C\cup \overline{D}\subset G.$

\medskip
We put $\widetilde{f}_n=\frac{1}{f_n(1/z)}.$ Since $f_n(z)=z+o(1)$
as $z\rightarrow\infty,$
$\lim\limits_{z\rightarrow\infty}f_n(z)=\infty.$ Set
$\widetilde{f}_n(0)=0.$ Set $\widetilde{f}_n(0):=\infty$. Since
$f_n(z)=z+o(1)$ as $z\rightarrow\infty,$ there is a neighborhood $U$
of the origin and a function $\varepsilon: U\rightarrow {\Bbb C}$
such that $f_n(1/z)=1/z+\varepsilon(1/z),$ where $z\in U$ and
$\varepsilon(1/z)\rightarrow 0$ as $z\rightarrow 0.$ Thus,
$$\frac{\widetilde{f}_n(\Delta z)-\widetilde{f}_n(0)}{\Delta z}=
\frac{1}{\Delta z}\cdot \frac{1}{1/(\Delta z)+\varepsilon(1/\Delta
z)}=\frac{1}{1+(\Delta z)\cdot \varepsilon(1/\Delta z)}\rightarrow 1
$$
as $\Delta z\rightarrow 0.$ This proves that there exists
$\widetilde{f}_n^{\,\prime}(0),$ and
$\widetilde{f}_n^{\,\prime}(0)=1.$ Since $\mu(z)$ vanishes outside
$G,$ the mapping $f_n$ is conformal in some neighborhood $V:={\Bbb
C}\setminus B(0, 1/r_0)$ of the infinity, and the number $1/r_0$
depends only on $G,$ and $G\subset B(0, 1/r_0).$ In this case, the
mapping $\widetilde{f}_n=\frac{1}{f_n(1/z)}$ is conformal in $B(0,
r_0).$ In addition, the mapping $F_n(z):=\frac{1}{r_0}\cdot
\widetilde{f}_n(r_0z)$ is a homeomorphism of the unit disk such that
$F_n(0)=0$ and $F_n^{\,\prime}(1)=1.$ By Koebe's theorem on~1/4
(see, e.g., \cite[Theorem~1.3]{CG}, cf.~\cite[Theorem~1.1.3]{GR})
$F_n({\Bbb D})\supset B (0, 1/4).$ Then
\begin{equation}\label{eq4C}
\widetilde{f}_n(B(0, r_0))\supset B(0, r_0/4)\,.
\end{equation}
By~(\ref{eq4C})
\begin{equation}\label{eq5B}
(1/f_n)(\overline{\Bbb C}\setminus \overline{B(0, 1/r_0)})\supset
B(0, r_0/4)\,.
\end{equation}
Taking into account formula~(\ref{eq5B}), we show that
\begin{equation}\label{eq6B}
f_n(\overline{\Bbb C}\setminus \overline{B(0, 1/r_0)})\supset
\overline{{\Bbb C}}\setminus \overline{B(0, 4/r_0)}\,.
\end{equation}
Indeed, let $y\in \overline{{\Bbb C}}\setminus \overline{B(0,
4/r_0)}.$ Now, $\frac{1}{y}\in B(0, r_0/4).$ By~(\ref{eq5B}),
$\frac{1}{y}=(1/f_n)(x),$ $x\in \overline{\Bbb C}\setminus
\overline{B(0, 1/r_0)}.$ Thus, $y=f_n(x),$ $x\in \overline{\Bbb
C}\setminus \overline{B(0, 1/r_0)},$ which proves~(\ref{eq6B}).

\medskip
Since $f_n$ is a homeomorphism in~${\Bbb C},$ by~(\ref{eq6B}) we
obtain that
\begin{equation}\label{eq2}
f_n(B(0, 1/r_0))\subset B(0, 4/r_0)\,.
\end{equation}
On the other hand, since $f_n$ are $n$-quasiconformal, the mappings
$g_n=f^{\,-1}_n$ are also quasiconformal; in particular, they belong
to the class~$W_{\rm loc}^{1, 2}(B(0, 1/r_0)).$
By~\cite[Theorem~6.10]{MRSY$_1$}, by~(\ref{eq10C}) and
by~(\ref{eq2})
\begin{equation} \label{eq2*B}
M(g_n(\Gamma))\leqslant \int\limits_{f_n(B(0,
1/r_0))}K_{\mu_{g_n}}(w)\cdot\rho_*^2 (w) \,dm(w)\leqslant
\int\limits_{B(0, 4/r_0)}Q(w)\cdot\rho_*^2 (w) \,dm(w)
\end{equation}
for any $n\in {\Bbb N},$ any path family $\Gamma$ in $f_n(B(0,
1/r_0)),$ and each function~$\rho_*\in {\rm adm}\,\Gamma.$ It
follows from~(\ref{eq2*B}) that
$$
M(g_n(\Gamma))\leqslant \int\limits_{B(0, 4/r_0)}Q(w)\cdot\rho_*^2
(w) \,dm(w)
$$
for the same functions $\rho.$ By Proposition~\ref{pr1} the family
$\{f_n\}_{n=1}^{\infty}$ is equicontinuous in~$B(0, 1/r_0).$ Thus,
by the Arzela-Ascoli theorem $f_n$ is a normal family of mappings
(see e.g.~\cite[Theorem~20.4]{Va}), in other words, there is a
subsequence $f_{n_k}$ of $f_n,$ converging locally uniformly in
${\Bbb C}$ to some map $f:B(0, 1/r_0)\rightarrow \overline{B(0,
4/r_0)}.$ Note also that $\mu_n(z)\rightarrow \mu(z)$ as
$n\rightarrow \infty$ and $\nu_n(z)\rightarrow \nu(z)$ as
$n\rightarrow \infty$ for almost all $z\in D,$ because $|\mu(z)|<1$
a.e. and, therefore, $K_{\mu}(z)$ in~(\ref{eq1D}) is finite for
almost all $z\in D.$ Then by~(\ref{eq10B}) and Lemma~\ref{pr2} the
map $f$ belongs to the class $W_{\rm loc}^{1, p}(D)$ and, in
addition, $f$ is a solution of~(\ref{eq2J}).

\medskip
Let us to prove that the limit mapping $f$ satisfies the condition
$f(z)=z+o(1)$ as $z\rightarrow\infty.$ Note that the family of
mappings
$F_{n_k}(z):=\frac{1}{r_0}\cdot\frac{1}{f_{n_k}(\frac{1}{r_0z})}$ is
compact in the unit disk (see, e.g., \cite[Theorem~1.10]{CG},
cf.~\cite[Theorem~1.2 Ch.~I]{GR}). Without loss of generality, we
may consider that $F_{n_k}$ converges locally uniformly in ${\Bbb
D}.$ Now $F(z)=\frac{1}{r_0}\cdot\frac{1}{f(\frac{1}{r_0z})}$
belongs to the class $S,$ consisting of conformal mappings $F$ of
the unit disk that satisfy the conditions $F(0)=0,$
$F^{\,\prime}(0)=1.$ Then the expansions of functions $F$ and
$F_{n_k}$ in a Taylor series at the origin have the form
\begin{equation}\label{eq1EA}F_{n_k}(z)=z+c_kz^2+z^2\cdot \varepsilon_k(z)\,,\quad
k=1,2,\ldots\,,
\end{equation}
\begin{equation}\label{eq2K}
F(z)=z+c_0z^2+z^2\cdot \varepsilon_0(z)\,,
\end{equation}
where $\varepsilon_k(z)$ and $\varepsilon_0(z)$ tend to zero as
$z\rightarrow 0.$ It follows from~(\ref{eq1EA}) and (\ref{eq2K})
that
\begin{equation}\label{eq3A}
f_{n_k}(t)=\frac{r_0t^2}{r_0t+c_k+\varepsilon_k\left(\frac{1}{r_0t}\right)}\,,\quad
f_{n_k}(t)-t=-\frac{c_k+\varepsilon_k\left(\frac{1}{r_0t}\right)}
{r_0+\frac{c_k}{t}+\frac{\varepsilon_k\left(\frac{1}{r_0t}\right)}{t}}\,,
\quad k=1,2,\ldots\,,
\end{equation}
\begin{equation}\label{eq4DA}
f(t)=\frac{r_0t^2}{r_0t+c_0+\varepsilon_0\left(\frac{1}{r_0t}\right)}\,,\quad
f(t)-t=-\frac{c_0+\varepsilon_0\left(\frac{1}{r_0t}\right)}
{r_0+\frac{c_0}{t}+\frac{\varepsilon_0\left(\frac{1}{r_0t}\right)}{t}}\,.
\end{equation}
In particular, passing to the limit in~(\ref{eq3A}) as
$t\rightarrow\infty,$ we obtain that $f_{n_k}(t)-t\rightarrow
-\frac{c_k}{r_0}.$ Since $f_{n_k}(t)-t\rightarrow 0$ as
$t\rightarrow \infty,$ we obtain that $c_k=0.$ By the Weierstrass
theorem on the convergence of the coefficients of the Taylor series
(see, e.g., \cite[Theorem~1.1.I]{Gol}) we obtain that
$c_k=0\rightarrow c_0$ as $k\rightarrow\infty.$ Thus, $c_0=0$
in~(\ref{eq4DA}), in other words, the mapping $f$ also has a
hydrodynamic normalization: $f(z)=z+o(1)$ as $z\rightarrow\infty.$
Theorem is proved.~$\Box$

\medskip
{\it Proof of Corollary~\ref{cor5}}. Let $G$ and $K$ be the same as
in the condition of the corollary.  Repeating the proof of
Theorem~\ref{th1A}, we observe that $g_n$ satisfy the relation
$$
M(g_n(\Gamma))\leqslant \int\limits_{f_n(B(0,
1/r_0))}K_{\mu_{g_n}}(w)\cdot\rho_*^2 (w) \,dm(w)\leqslant
\int\limits_{\Bbb C}Q^{\,\prime}(w)\cdot\rho_*^2 (w) \,dm(w)\,,
$$
where $Q^{\,\prime}(w)=\begin{cases}Q(w)\,,& w\in B(0, 4/r_0)\,,\\
0\,,& \text{otherwise}\end{cases},$ where the function
$Q^{\,\prime}$ is integrable in ${\Bbb C}.$ Now,
by~\cite[Theorem~4.1]{SevSkv}
the relation
$$
|f_n(x)-f_n(y)|\leqslant\frac{C}{\log^{1/2}\left(1+\frac{r_0}{2|x-y|}\right)}
$$
holds for any $x, y\in K$ and any compact set $K$ in $G,$ where
$C=C(K, \Vert Q^{\,\prime}\Vert_1, G)>0$ is some constant depending
only on $K,$ $G$ and $\Vert Q^{\,\prime}\Vert_1,$ $\Vert
Q^{\,\prime}\Vert_1$ denotes $L^1$-norm of $Q^{\,\prime}$ in ${\Bbb
R}^n,$ and $r_0=d(K,
\partial G).$ Passing here to the limit as
$n\rightarrow\infty,$ we obtain the required
relation~(\ref{eq2E}).~$\Box$

\section{Convergence theorem}

Recall that, a mapping $f:D\rightarrow {\Bbb R}^n,$ $n\geqslant 2,$
is called a {\it mapping with a finite distortion,} if $f\in W_{\rm
loc}^{1, 1}(D)$ and there is a function $K:D\rightarrow [1, \infty)$
such that $\Vert f^{\,\prime}(x)\Vert^{n}\leqslant K(x)\cdot |J(x,
f)|$ for a.a. $x\in D.$ Observe that, if $f$ has a finite
distortion, then $g:=f^{\,-1}$ is of finite distortion, as well (see
\cite[Theorem~1.2]{HK}).

\medskip
Given a function $Q:{\Bbb C}\rightarrow {\Bbb C},$ numbers
$1<p\leqslant 2,$ $0<M<\infty$ and a domain $G$ such that
$\overline{G}$ is a compactum in ${\Bbb C},$ denote by $\frak{F}_{Q,
p, M}(G)$ a family of all $W_{\rm loc}^{1, 1}$-homeomorphisms
$f:{\Bbb C}\rightarrow {\Bbb C}$ with a finite distortion which are
conformal outside $G$ and are solutions of the equation~(\ref{eq2J})
in ${\Bbb C}$ such that

\medskip
1) $f(z)=z+\varepsilon_f(z),$ where $\varepsilon_f(z)\rightarrow 0$
as $z\rightarrow\infty ;$

\medskip
2) there exist a number $1<p\leqslant 2$ such that, for any bonded
domain $G\subset {\Bbb C}$ there exists a constant $M=M_G>0$ such
that
\begin{equation}\label{eq10D}
\int\limits_{G}K_{I, p}(w, g)\,dm(w)\leqslant M
\end{equation}
for all $n=1,2,\ldots ,$ where $g:=f^{\,-1}$ and $K_{I, p}(w, g)$ is
defined in~(\ref{eq18});

\medskip
2) the inequality
\begin{equation}\label{eq10E}
K_{\mu_g}(w)\leqslant Q(w)
\end{equation}
holds for $g:=f^{\,-1}$ and for a.e. $w\in {\Bbb C},$ where
$K_{\mu_{g}}$ is defined in~(\ref{eq3D}).

\medskip
The following statement holds.

\medskip
\begin{theorem}\label{th2}
{\sl\, Assume that, for each $0<r_1<r_2<1$ and $y_0\in {\Bbb C}$
there is a set $E\subset[r_1, r_2]$ of positive linear Lebesgue
measure such that the function $Q$ is integrable over the circles
$S(y_0, r)$ for any $r\in E.$ Then $\frak{F}_{Q, p, M}(G)$ is a
normal family. If $f_n\in \frak{F}_{Q, p, M}(G),$ $n=1,2,\ldots ,$
$f_n\rightarrow f$ as $n\rightarrow\infty$ locally uniformly, and
$\mu_n\rightarrow \mu$ and $\nu_n\rightarrow\nu$ as
$n\rightarrow\infty,$ then $f$ satisfies the equation~(\ref{eq2J}).
In this case, $f(z)=z+\varepsilon(z),$ $\varepsilon(z)\rightarrow 0$
as $z\rightarrow\infty,$ as well.}
\end{theorem}

\medskip
\begin{corollary}\label{cor6}
{\sl\, In particular, the statement of Theorem~\ref{th2} holds, if
instead of the conditions on the function $Q,$ specified in this
theorem, we require that $Q\in L_{\rm loc}^1.$}
\end{corollary}

\medskip
{\it Proof of Theorem~\ref{th2}.} Note that the proof follows the
same principle as the proof of the previous Theorem~\ref{th1A}.

Let us to prove that $\frak{F}_{Q, p, M}(G)$ forms a normal family
of mappings. Fix an arbitrary domain $G\subset {\Bbb C}$ with a
compct closure.

\medskip
Given a mapping $f\in\frak{F}_{Q, p, M}(G),$ we put
$\widetilde{f}=\frac{1}{f(1/z)}.$ Since $f(z)=z+\varepsilon(z),$
$\varepsilon(z)\rightarrow 0$ as $z\rightarrow\infty,$ we obtain
that $\lim\limits_{z\rightarrow\infty}f(z)=\infty.$ Set
$\widetilde{f}(0)=0.$ Set $\widetilde{f}(0):=\infty$. Now, there is
a neighborhood $U$ of the origin such that
$f(1/z)=1/z+\varepsilon(1/z),$ where $z\in U$ and
$\varepsilon(1/z)\rightarrow 0$ as $z\rightarrow 0.$ Thus,
$$\frac{\widetilde{f}(\Delta z)-\widetilde{f}(0)}{\Delta z}=
\frac{1}{\Delta z}\cdot \frac{1}{1/(\Delta z)+\varepsilon(1/\Delta
z)}=\frac{1}{1+(\Delta z)\cdot \varepsilon(1/\Delta z)}\rightarrow 1
$$
as $\Delta z\rightarrow 0.$ This proves that there exists
$\widetilde{f}^{\,\prime}(0),$ and $\widetilde{f}^{\,\prime}(0)=1.$
Since $\mu(z)$ vanishes outside $G,$ the mapping $f$ is conformal in
some neighborhood $V:={\Bbb C}\setminus B(0, 1/r_0)$ of the
infinity, and the number $1/r_0$ depends only on $G,$ and $G\subset
B(0, 1/r_0).$ In this case, the mapping
$\widetilde{f}=\frac{1}{f(1/z)}$ is conformal in $B(0, r_0).$ In
addition, the mapping $F(z):=\frac{1}{r_0}\cdot \widetilde{f}(r_0z)$
is a homeomorphism of the unit disk such that $F(0)=0$ and
$F^{\,\prime}(1)=1.$ By Koebe's theorem on~1/4 (see, e.g.,
\cite[Theorem~1.3]{CG}, cf.~\cite[Theorem~1.1.3]{GR}) $F({\Bbb
D})\supset B (0, 1/4).$ Then
\begin{equation}\label{eq4E}
\widetilde{f}(B(0, r_0))\supset B(0, r_0/4)\,.
\end{equation}
By~(\ref{eq4E})
\begin{equation}\label{eq5C}
(1/f)(\overline{\Bbb C}\setminus \overline{B(0, 1/r_0)})\supset B(0,
r_0/4)\,.
\end{equation}
Taking into account formula~(\ref{eq5C}), we show that
\begin{equation}\label{eq6C}
f(\overline{\Bbb C}\setminus \overline{B(0, 1/r_0)})\supset
\overline{{\Bbb C}}\setminus \overline{B(0, 4/r_0)}\,.
\end{equation}
Indeed, let $y\in \overline{{\Bbb C}}\setminus \overline{B(0,
4/r_0)}.$ Now, $\frac{1}{y}\in B(0, r_0/4).$ By~(\ref{eq5C}),
$\frac{1}{y}=(1/f)(x),$ $x\in \overline{\Bbb C}\setminus
\overline{B(0, 1/r_0)}.$ Thus, $y=f(x),$ $x\in \overline{\Bbb
C}\setminus \overline{B(0, 1/r_0)},$ which proves~(\ref{eq6C}).

\medskip
Since $f$ is a homeomorphism in~${\Bbb C},$ by~(\ref{eq6C}) we
obtain that
\begin{equation}\label{eq2A}
f(B(0, 1/r_0))\subset B(0, 4/r_0)\,.
\end{equation}

\medskip
On the other hand, taking into account the comments made above, each
mapping $g=f^{\,-1},$ has a finite distortion. Then
by~\cite[Lemma~3.1 and Proposition~2.1]{LSS}
$$M(g(\Sigma_{r_1, r_2}))\geqslant \int\limits_{r_1}^{r_2}
\frac{dr}{\Vert Q\Vert_1(r)}\,,$$
where $\Sigma_{r_1, r_2}$ denotes the family of all circles $S(y_0,
r)\cap f(G),$ $r\in (r_1, r_2),$ and $\Vert
Q\Vert_1(r)=\int\limits_{S(y_0, r)\cap
f(G)}Q(y)\,\,d\mathcal{H}^{1}(y).$ Applying the Ziemer and Hesse
theorems, see~\cite[Theorem~3.13]{Zi} and~\cite[Theorem~5.5]{Hes},
we obtain that
$$M(g(\Gamma(S(y_0, r_1), S(y_0, r_2), f(G))))\leqslant \frac{2\pi}{\int\limits_{r_1}^{r_2}
\frac{dr}{rq_{y_0}(r)}}\,,$$
or, in some another form,
$$M(\Gamma_f(y_0, r_1, r_2))\leqslant \frac{2\pi}{\int\limits_{r_1}^{r_2}
\frac{dr}{rq_{y_0}(r)}}\,.$$
Now, by Remark~\ref{rem1} it follows that $\frak{F}_{Q, p, M}(G)$ is
equicontinuous in $G.$ Finally, by the Arzela-Ascoli theorem
$\frak{F}_{Q, p, M}(G)$ is a normal family of mappings (see
e.g.~\cite[Theorem~20.4]{Va}).

\medskip
Assume now that $f_n\in \frak{F}_{Q, p, M}(G),$ $n=1,2,\ldots ,$ is
a sequence converging locally uniformly in ${\Bbb C}$ to some map
$f:{\Bbb C}\rightarrow {\Bbb C}.$ Note also that
$\mu_n(z)\rightarrow \mu(z)$ as $n\rightarrow \infty$ and
$\nu_n(z)\rightarrow \nu(z)$ as $n\rightarrow \infty$ for almost all
$z\in D,$ because $|\mu(z)|<1$ a.e. and, therefore, $K_{\mu}(z)$
in~(\ref{eq1D}) is finite for almost all $z\in D.$ Then
by~(\ref{eq10B}) and Lemma~\ref{pr2} the map $f$ belongs to the
class $W_{\rm loc}^{1, p}(D)$ and, in addition, $f$ is a solution
of~(\ref{eq2J}).

\medskip
Let us to prove that the limit mapping $f$ satisfies the condition
$f(z)=z+o(1)$ as $z\rightarrow\infty.$ Note that the family of
mappings
$F_{n_k}(z):=\frac{1}{r_0}\cdot\frac{1}{f_{n_k}(\frac{1}{r_0z})}$ is
compact in the unit disk (see, e.g., \cite[Theorem~1.10]{CG},
cf.~\cite[Theorem~1.2 Ch.~I]{GR}). Without loss of generality, we
may consider that $F_{n_k}$ converges locally uniformly in ${\Bbb
D}.$ Now $F(z)=\frac{1}{r_0}\cdot\frac{1}{f(\frac{1}{r_0z})}$
belongs to the class $S,$ consisting of conformal mappings $F$ of
the unit disk that satisfy the conditions $F(0)=0,$
$F^{\,\prime}(0)=1.$ Then the expansions of functions $F$ and
$F_{n_k}$ in a Taylor series at the origin have the form
\begin{equation}\label{eq1EB}F_{n_k}(z)=z+c_kz^2+z^2\cdot \varepsilon_k(z)\,,\quad
k=1,2,\ldots\,,
\end{equation}
\begin{equation}\label{eq2KB}
F(z)=z+c_0z^2+z^2\cdot \varepsilon_0(z)\,,
\end{equation}
where $\varepsilon_k(z)$ and $\varepsilon_0(z)$ tend to zero as
$z\rightarrow 0.$ It follows from~(\ref{eq1EB}) and (\ref{eq2KB})
that
\begin{equation}\label{eq3B}
f_{n_k}(t)=\frac{r_0t^2}{r_0t+c_k+\varepsilon_k\left(\frac{1}{r_0t}\right)}\,,\quad
f_{n_k}(t)-t=-\frac{c_k+\varepsilon_k\left(\frac{1}{r_0t}\right)}
{r_0+\frac{c_k}{t}+\frac{\varepsilon_k\left(\frac{1}{r_0t}\right)}{t}}\,,
\quad k=1,2,\ldots\,,
\end{equation}
\begin{equation}\label{eq4F}
f(t)=\frac{r_0t^2}{r_0t+c_0+\varepsilon_0\left(\frac{1}{r_0t}\right)}\,,\quad
f(t)-t=-\frac{c_0+\varepsilon_0\left(\frac{1}{r_0t}\right)}
{r_0+\frac{c_0}{t}+\frac{\varepsilon_0\left(\frac{1}{r_0t}\right)}{t}}\,.
\end{equation}
In particular, passing to the limit in~(\ref{eq3B}) as
$t\rightarrow\infty,$ we obtain that $f_{n_k}(t)-t\rightarrow
-\frac{c_k}{r_0}.$ Since $f_{n_k}(t)-t\rightarrow 0$ as
$t\rightarrow \infty,$ we obtain that $c_k=0.$ By the Weierstrass
theorem on the convergence of the coefficients of the Taylor series
(see, e.g., \cite[Theorem~1.1.I]{Gol}) we obtain that
$c_k=0\rightarrow c_0$ as $k\rightarrow\infty.$ Thus, $c_0=0$
in~(\ref{eq4F}), in other words, the mapping $f$ also has a
hydrodynamic normalization: $f(z)=z+o(1)$ as $z\rightarrow\infty.$
Theorem is proved.~$\Box$

\medskip
\begin{example}\label{ex1}
Let $p=2,$ let $q\geqslant 1$ be an arbitrary number and let
$0<\alpha<2/q.$ As usual, we use the notation $z=re^{i\theta},$
$r\geqslant 0,$ $\theta\in [0, 2\pi).$ Put
\begin{equation}\label{eq3A*}\mu(z)= \left
\{\begin{array}{rr}
 e^{2i\theta}\frac{2r-\alpha(2r-1)}{2r+\alpha(2r-1)},& 1/2<|z|<1,
\\ 0\ , & |z|\leqslant 1/2\,.
\end{array} \right.
\end{equation}
Using the ratio
$$\mu_f(z)=\frac{\overline{\partial} f}{\partial f}=e^{2i\theta}\frac{rf_r+if_{\theta}}{rf_r-if_{\theta}}\,,$$
see~(11.129) in \cite{MRSY$_2$}, we obtain that the mapping
\begin{equation}\label{eq4A*}f(z)=\left
\{\begin{array}{rr}
 \frac{z}{|z|}(2|z|-1)^{1/\alpha},& 1/2<|z|<1,
\\ 0\ , & |z|\leqslant 1/2
\end{array} \right.
\end{equation}
is a solution of the equation~$f_{\overline{z}}=\mu(z)\cdot f_z,$
where $\mu$ is defined by~(\ref{eq3A*}). Note that for $\mu$
in~(\ref{eq3A*}), the corresponding maximal dilatation $K_{\mu}$ is
the function
\begin{equation}\label{eq5A}K_{\mu}(z)=\left
\{\begin{array}{rr}
 \frac{2|z|}{\alpha(2|z|-1)},& 1/2<|z|<1,
\\ 1\ , & |z|\leqslant 1/2
\end{array} \right.\,.
\end{equation}
Let $k>1/\alpha.$ Observe that $K_{\mu}(z)\leqslant k$ for
$|z|\geqslant \frac{1}{2}\cdot \frac{k\alpha}{k\alpha-1}$ and
$K_{\mu}(z)>k$ otherwise. As above, we set
\begin{equation*}\label{eq5E}\mu_k(z)= \left
\{\begin{array}{rr}
 \mu(z),& K_{\mu}(z)\leqslant k,
\\ 0\ , & K_{\mu}(z)> k\,.
\end{array} \right.
\end{equation*}
Observe that the mappings
\begin{equation*}\label{eq6CB}f_k(z)=\left
\{\begin{array}{rr}
 \frac{z}{|z|}(2|z|-1)^{1/\alpha},& \frac{1}{2}\cdot \frac{k\alpha}{k\alpha-1}<|z|<1,
\\ \frac{z}{\frac{1}{2}\left(\frac{k\alpha}{k\alpha-1}\right)
}\cdot{\left(\frac{1}{k\alpha-1}\right)}^{1/\alpha}\ , &
|z|\leqslant \frac{k\alpha}{k\alpha-1}
\end{array} \right.\,,
\end{equation*}
are homeomorphic solutions of the equation
$f_{\overline{z}}=\mu_k(z)\cdot f_z.$ Besides that, the inverse
mappings $g_k(y)=f_k^{\,-1}(y)$ are calculated by the relations
\begin{equation}\label{eq7E}g_k(y)=\left
\{\begin{array}{rr}
 \frac{y(|y|^{\alpha}+1)}{2|y|},& \left(\frac{k\alpha}{k\alpha-1}-1\right)^{1/\alpha}<|y|<1,
\\ \frac{y\cdot\frac{k\alpha}{2(k\alpha-1)}}
{\left(\frac{k\alpha}{k\alpha-1}-1\right)^{1/\alpha}} , &
|y|\leqslant\left(\frac{k\alpha}{k\alpha-1}-1\right)^{1/\alpha}
\end{array} \right.\,.
\end{equation}
It follows from~(\ref{eq5A}) that
\begin{equation}\label{eq7F}K_{\mu_k}(z)=\left
\{\begin{array}{rr}
 \frac{4|z|}{2\alpha(2|z|-1)},& \frac{1}{2}\cdot \frac{k\alpha}{k\alpha-1}<|z|<1,
\\ 1\ , & |z|\leqslant \frac 12\cdot\frac{k\alpha}{k\alpha-1}
\end{array} \right.\,.
\end{equation}
We should check that relation~(\ref{eq10B}) holds for some function
$Q$ that is integrable in ${\Bbb D}.$ For this purpose, we
substitute the maps $g_k$ from~(\ref{eq7E}) into the maximal
dilatation~$K_{\mu_k}$ defined by the equality~(\ref{eq7F}). Then
\begin{equation*}\label{eq8C}K_{\mu_{g_k}}(y)=\left
\{\begin{array}{rr}
 \frac{|y|^{\alpha}+1}{\alpha|y|^{\alpha}},&
 \left(\frac{k\alpha}{k\alpha-1}-1\right)^{1/\alpha}<|y|<1\,,
\\ 1\ , & |y|\leqslant\left(\frac{k\alpha}{k\alpha-1}-1\right)^{1/\alpha}
\end{array} \right.\,.
\end{equation*}
Note that $K_{\mu_{g_k}}(y)\leqslant Q(y):=
\frac{|y|^{\alpha}+1}{\alpha|y|^{\alpha}}$ for all $y\in {\Bbb D}.$
Moreover, the function $Q$ is integrable in ${\Bbb D}$ even in the
degree $q,$ and not only in the degree 1 (see the arguments used in
considering~\cite[Proposition~6.3]{MRSY$_2$}). We extend each of the
mappings $f_k$ identically to the whole plane, and set $Q(y)\equiv
1$ for $y\not\in {\Bbb C}.$ It follows from the considerations
mentioned above that, all mappings $f_k,$ $k=1,2,\ldots,$ belong to
some class $\frak{F}_{Q, p, M}(G)$ with $Q$ mentioned above. Note
that the family $f_k,$ $k=1,2,\ldots,$ is normal, but not compact,
since the limit mapping $f$ of this sequence is not a homeomorphism.
In order for the compactness of the class to be satisfied, other
conditions on the characteristics of mappings are necessary, which
will be considered below.
\end{example}

\medskip
Let $D$ be a domain in ${\Bbb C}.$ Suppose that a function
${\varphi}:D\rightarrow{\Bbb R}$ is locally integrable in some
neighborhood of a point $z_0\in D .$ We say that ${\varphi}$ has a
{\it finite mean oscillation} at $z_0\in D,$ and we write
$\varphi\in FMO(z_0),$ if the relation
\begin{equation*}\label{eq17:}
{\limsup\limits_{\varepsilon\rightarrow
0}}\,\frac{1}{\pi\varepsilon^2}\int\limits_{B(z_0,\,\varepsilon)}
|{\varphi}(z)-\overline{{\varphi}}_{\varepsilon}|\ dm(z)\, <\,
\infty
\end{equation*}
holds, where
$\overline{{\varphi}}_{\varepsilon}=\frac{1}{\pi\varepsilon^2}\int\limits_{B(
z_0,\,\varepsilon)} {\varphi}(z)\ dm(z)$ (see, e.g.,
\cite[section~2]{RSY$_1$}).
We say that a function $\varphi$  has a finite mean oscillation in
$D,$ and we write $\varphi\in FMO(D),$ if $\varphi\in FMO(z_0)$ for
any $z_0\in D.$

\medskip
Given a function $Q:{\Bbb C}\rightarrow {\Bbb C}$ and a domain $G$
with a compact closure in ${\Bbb C},$ denote by $\frak{R}_{Q}(G)$ a
family of all $W_{\rm loc}^{1, 1}$-homeomorphisms $f:{\Bbb
C}\rightarrow {\Bbb C}$ with a finite distortion which are conformal
outside $G$ and are solutions of the equation~(\ref{eq2J}) in ${\Bbb
C}$ such that

\medskip
1) $f(z)=z+\varepsilon_f(z),$ where $\varepsilon_f(z)\rightarrow 0$
as $z\rightarrow\infty ;$

\medskip
2) the inequality
\begin{equation}\label{eq10F}
K_{\mu_g}(w)\leqslant Q(w)
\end{equation}
holds for $g:=f^{\,-1}$ and for a.e. $w\in {\Bbb C},$ where
$K_{\mu_{g}}$ is defined in~(\ref{eq3D}).

The following statement holds.

\medskip
\begin{theorem}\label{th1}
{\sl If 1) $Q(z)\in FMO({\Bbb C}),$ or 2) $Q\in L_{\rm loc}^1({\Bbb
C})$ and, in addition,
\begin{equation}\label{eq5**}
\int\limits_{0}^{\delta(w_0)}\frac{dt}{tq_{w_0}(t)}=\infty
\end{equation}
for any~$w_0\in {\Bbb C}$ and some~$\delta(w_0)>0,$
$q_{w_0}(r)=\frac{1}{2\pi}\int\limits_{0}^{2\pi}Q(w_0+re^{\,i\theta})\,d\theta,$
then the class $\frak{R}_{Q}(G)$ is compact.}
\end{theorem}

\medskip
\begin{proof}
Since the normality of $\frak{R}_{Q}(G)$ follows by
Theorem~\ref{th2}, we need to prove the closeness of
$\frak{R}_{Q}(G).$ Let $f_n\in \frak{R}_{Q}(G),$ $n=1,2,\ldots,$ and
let $f_n\rightarrow f$ as $n\rightarrow \infty.$ Since
$f_n(z)=z+o(1)$ as $z\rightarrow\infty,$ denoting by
$g_n:=f_n^{\,-1}$ and $f_n(z)=w$ we obtain that
$z=g_n(w)=w-o(1)=w-\varepsilon(z(w))=w+\varepsilon_1(w),$ where
$\varepsilon_1(w)\rightarrow 0$ as $w\rightarrow \infty.$ Arguing
similarly to the proof of Theorem~\ref{th1A}, we obtain that
$g_n(B(0, 1/r_0))\subset B(0, 4/r_0)$ for any $n\in {\Bbb N}$ and
any sufficiently small $r_0>0.$ Then the sequence $g_n:=f_n^{\,-1}$
forms an equicontinuous family of mappings, as well
(see~\cite[Theorems~6.1 and 6.5]{RS}). Therefore, by the
Arzela-Ascoli theorem $g_k$ is a normal family (see
e.g.~\cite[Theorem~20.4]{Va}), in other words, there is a
subsequence $g_{n_k}$ of $g_n$ converging locally uniformly in
${\Bbb C}$ to some map $g:{\Bbb C}\rightarrow \overline{{\Bbb C}}.$
Arguing as in the proof of Theorem~\ref{th1A}, we obtain that $g(z)$
has a hydrodynamical normalization at the infinity.  Then, by virtue
of~\cite[Theorems~4.1, 4.2]{RSS} the mapping $g$ is a homeomorphism
in ${\Bbb C}$ and $g:{\Bbb C}\rightarrow {\Bbb C}.$ Moreover,
$g({\Bbb C})={\Bbb C}$ because $g$ is a homeomorphism and the point
$\infty$ is isolated boundary point of $g$ (see, e.g.,
\cite[Theorem~6.2]{MRSY$_2$}). Thus, by~\cite[Lemma~3.1]{RSS}, we
also have that $g_{n_k}\rightarrow f=g^{\,- 1}$ as
$k\rightarrow\infty$ locally uniformly in ${\Bbb C}.$ Observe that,
since $f_{n_k}$ has a finite distortion, then
$g_{n_k}:=f^{\,-1}_{n_k}$ is of finite distortion, as well (see
\cite[Theorem~1.2]{HK}). Now, $K_{\mu_g}(w)\leqslant Q(w)$ for
almost any $w\in {\Bbb C}$ by \cite[Theorem~16.1]{RSS}.~$\Box$
\end{proof}


\medskip
{\bf \noindent Oleksandr Dovhopiatyi} \\
{\bf 1.} Zhytomyr Ivan Franko State University,  \\
40 Bol'shaya Berdichevskaya Str., 10 008  Zhytomyr, UKRAINE \\
alexdov1111111@gmail.com

\medskip
{\bf \noindent Evgeny Sevost'yanov} \\
{\bf 1.} Zhytomyr Ivan Franko State University,  \\
40 Bol'shaya Berdichevskaya Str., 10 008  Zhytomyr, UKRAINE \\
{\bf 2.} Institute of Applied Mathematics and Mechanics\\
of NAS of Ukraine, \\
1 Dobrovol'skogo Str., 84 100 Slavyansk,  UKRAINE\\
esevostyanov2009@gmail.com

\end{document}